\documentclass[11pt]{article}

\usepackage{amsmath}    
\usepackage{graphicx}   
\usepackage{verbatim}   
\usepackage{color}      
\usepackage{subfigure}  
\usepackage{amssymb,tikz}
\usepackage{soul}

\newtheorem{theorem}{Theorem}[section]

\newtheorem{conjecture}[theorem]{Conjecture}

\newtheorem{lemma}[theorem]{Lemma}

\newcommand{\proof}{\noindent{\bf Proof.\ }}
\newcommand{\qed}{\hfill $\square$\medskip}

\begin{document}

\title{Revisiting path-type covering and partitioning problems}

\author{
	Paul Manuel
}
{\Large This is a survey article which is at the initial stage. The author will appreciate to receive your comments and contributions to improve the quality of the article. The author's contact address is pauldmanuel@gmail.com.}
\date{\today}

\maketitle
\begin{center}

	Department of Information Science, \\
	College of Computing Science and Engineering, \\
	Kuwait University, Kuwait \\
	{\tt pauldmanuel@gmail.com}\\
	\medskip
		
\end{center}

\begin{abstract}
Covering problems belong to the foundation of graph theory. There are several types of covering problems in graph theory such as covering the vertex set by stars (domination problem), covering the vertex set by cliques (clique covering problem), covering the vertex set by independent sets (coloring problem), and covering the vertex set by paths or cycles. A similar concept which is partitioning problem is also equally important. Lately research in graph theory has produced unprecedented growth because of its various application in engineering and science. The covering and partitioning problem by paths itself have produced a sizable volume of literatures. The research on these problems is expanding in multiple directions and the volume of research papers is exploding. It is the time to simplify and unify the literature on different types of the covering and partitioning problems. The problems considered in this article are path cover problem, induced path cover problem, isometric path cover problem, path partition problem, induced path partition problem and isometric path partition problem. The objective of this article is to summarize the recent developments on these problems, classify their literatures and correlate the inter-relationship among the related concepts.
\end{abstract}


\medskip
\noindent{\bf AMS Subj.\ Class.: 05C12, 05C70, 05C76; 68Q17} 

\section{Introduction}
An undirected connected graph is represented by $G(V,E)$ where $V$ is the vertex set and $E$ is the edge set. There are several types of covering problems such as covering the vertex set by stars (domination problem), covering the vertex set by cliques (the clique covering problem), covering the vertex set by independent sets (the coloring problem), and covering the vertex set by path-types. A \emph{path-type} means a path, an induced path or an isometric path. A path $P$ is an \emph{induced path} in $G$ if the subgraph induced by the vertices of $P$ is a path. An induced path is also called a \emph{chordless path} ~\cite{Mezz10}. Let us recall that \emph{isometric path} and \emph{geodesic} are other names for shortest path. 

Mathematicians agree that the term  ``covering'' is over-used or inappropriately-used in graph theory. In set theory, there are two different terminologies: set covering problem and set hitting problem. The edge cover problem is a special case of set covering problem while the vertex cover problem is a special case of set hitting problem. While set theorists distinguished the two concepts by ``covering'' and ``hitting'', graph theorists chose to use the same term ``covering'' to mean two different concepts. Thus, right in the beginning, the term ``cover'' was over-used. The trend has been continuing. It is not uncommon that researchers use the term ``path cover problem'' to mean ``path partition problem'' and use the same term to mean ``induced path cover problem''. It is observed that the concepts ``cover'' \& ``partition'' and the terms ``path'' \& ``induced path'' are used inappropriately. There is no consistency in the usage of notations too. For example, while some authors use the notation $\pi$ to represent the path cover number, others use the notation $\mu$ to represent the same concept of path cover number. These ``inconsistent notation'' and ``inappropriate terminologies'' lead novice readers to confusion. This situation demands an article that classifies and simplifies the notations and terminologies of path-type covering and partitioning problems.  

In this article, we divide the ``theory of decomposition'' into six problems and classify the literatures accordingly. We also bring their related concepts together under one umbrella. This article targets young researchers in graph theory. The benefit of this article is that it exposes the gap and open space in research area under each problem. The article highlights that there is plenty of research to be done in this topic. In addition, this article will save time for young researchers to choose a right research topic in this domain and it may be a reference to apply common notations and appropriate terminologies in their research articles. 

In section \ref{sec:some_gen_concepts}, some relevant problems such as the $k$-path cover problem, the $S$-path cover problem and $k$-path vertex cover problem are discussed. Since these topics are not our main focus, we do not provide elaborate literature survey for these problems. In section \ref{sec:dis_Gallai_Milgram_thm}, the significance of Gallai-Milgram theorem and Berge path partition conjecture is highlighted. The three cover problems are discussed in section \ref{sec:three_cover_prob} and the three partition problems are studied in section \ref{sec:three_part_prob}.
\section{Path-type covering and partitioning problems}
\label{sec:path-type-cover-partition}
\begin{figure}[ht!]
	\begin{center}
		\scalebox{0.53}{\includegraphics{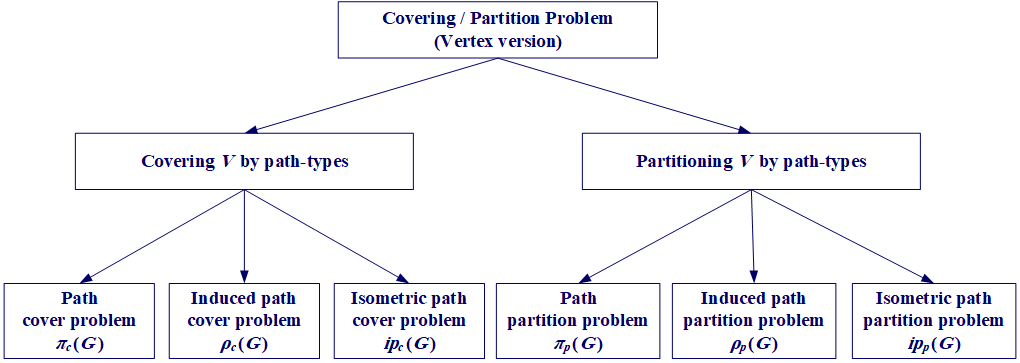}}
	\end{center}
	\caption{Six types of covering / partitioning problems}
	\label{F6TyProb}
\end{figure}

Brause and Krivo\v{s}-Bellu\v{s} \cite{BrKr17} state that the cover problem and the partition problem are central problems in graph theory. One such well-known concept in graph theory is ``clique covering'' and ``clique partitioning'' \cite{ErFa88}. A \emph{path cover} of $G$ is a set of paths such that every vertex $v \in V$ belongs to at least one path, whereas a \emph{path partition} of $G$ is a set of paths such that every vertex $v \in V$ belongs to exactly one path \cite{Hart06}. In a path partition, the paths are mutually vertex-disjoint. A few authors use the term ``vertex-disjoint path cover'' to mean ``path partition'' \cite{BaFa04}. The following six invariants (See Figure  \ref{F6TyProb}) are defined here:
\begin{enumerate}
	\item[(1)] The \emph{path cover number} of $G$, denoted by  $\pi_{c}(G)$, is the cardinality of a minimum path cover of $G$. The \emph{path cover problem} is to find a path cover of minimum cardinality in $G$. 
	\item[(2)] The \emph{induced path cover number} of $G$, denoted by  $\rho_{c}(G)$, is the cardinality of a minimum induced path cover of $G$. The \emph{induced path cover problem} is to find an induced path cover of minimum cardinality in $G$.
	\item[(3)] The \emph{isometric path cover number} of $G$, denoted by $ip_{c}(G)$, is the cardinality of a minimum isometric path cover of $G$. The \emph{isometric path cover problem} is to find an isometric path cover of minimum cardinality in $G$.
	\item[(4)] The \emph{path partition number} of $G$, denoted by $\pi_{p}(G)$, is the cardinality of a minimum path partition of $G$. The \emph{path partition problem} is to find a path partition of minimum cardinality in $G$.
	\item[(5)] The \emph{induced path partition number} of $G$, denoted by $\rho_{p}(G)$, is the cardinality of a minimum induced path partition of $G$. The \emph{induced path partition problem} is to find an induced path partition of minimum cardinality in $G$.
	\item[(6)] The \emph{isometric path partition number} of $G$, denoted by $ip_{p}(G)$, is the cardinality of a minimum isometric path partition of $G$. The \emph{isometric path partition problem} is to find an isometric path partition of minimum cardinality in $G$.
\end{enumerate}

While some authors \cite{ArRa90, Char94, ErFa88, Hart88, Hart06, PaCh07, Skup74, SriSu93, Wong99, YaCh94} use the term ``path partition'', others \cite{AsNi07a, BaFa04, CaCe12, CoDa13, HuCh06, MaMa09, NaOl03, Ore61, Row11, TaFa14, Yu17} use the term ``path cover'' (or vertex-disjoint path cover) in the place of ``path partition''. Some authors use the term ``decomposition'' for ``partition'' \cite{ArHa13, RaRa13}. 

A simple example in Figure \ref{FSimpEx} demonstrates the different behaviors of six different problems:
\begin{figure}[ht!]
	\begin{center}
		\scalebox{0.8}{\includegraphics{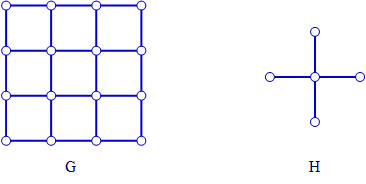}}
	\end{center}
	\caption{Illustration how the six problems are different}
	\label{FSimpEx}
\end{figure}  

$\pi_{c}(G) = 1$;
$\rho_{c}(G) = 2$; 
$ip_{c}(G) = 3$; 		
$\pi_{p}(G) = 1$; 
$\rho_{p}(G) = 2$; 
$ip_{p}(G) = 4$;
 
$\pi_{c}(H) = 2$; 
$\rho_{c}(H) = 2$; 
$ip_{c}(H) = 2$; 		
$\pi_{p}(H) = 3$; 
$\rho_{p}(H) = 3$; 
$ip_{p}(H) = 3$;

\noindent
The inequalities in Lemma \ref{L6Inq} are true because of the following facts:
\begin{enumerate}
	\item[(1)] An isometric path is an induced path. An induced path is a path.
	\item[(2)] A path partition is a path cover.
\end{enumerate}

\begin{lemma}
\label{L6Inq}
	For a given graph $G$,
	\begin{enumerate}
		\item[(1)] $\pi_{c}(G) \leq \rho_{c}(G) \leq ip_{c}(G)$; 	
		$\pi_{p}(G) \leq \rho_{p}(G) \leq ip_{p}(G)$;
		
		\item[(2)]	$\pi_{c}(G)\leq \pi_{p}(G)$; 
		$\rho_{c}(G)\leq \rho_{p}(G)$; 
		$ip_{c}(G)\leq ip_{p}(G)$; 
	\end{enumerate}
\end{lemma}

\section{Some generalized concepts}
\label{sec:some_gen_concepts}
In this section, we discuss some problems which are closely related to the six problems.
\subsection{The $k$-path cover problem}
\label{subsec:k-path_cover_prob}
A \emph{$k$-path} is a path having at most $k$ vertices. Here is a word of caution: While defining a $k$-path, some authors define $k$-path as a path of at most $k$ vertices \cite{JiLi06, Stei03} and some authors define $k$-path as a path of length at most $k$ \cite{Korp18}.  Korpelainen \cite{Korp18} states that the $3$-path partition problem is NP-complete on comparability graphs. According to his definition of $k$-path, the correct statement is that the $2$-path partition problem is NP-complete on comparability graphs.

A \emph{$k$-path cover} is a set of $k$-paths that cover $V$. The \emph{$k$-path cover problem} is to find a $k$-path cover  of minimum cardinality in $G$. This is a generalized version of path cover problem \cite{Stei03}. When $k = 2$, the $k$-path cover problem becomes the edge cover problem. It is straightforward to define other five problems in the same way \cite{Korp18, MoTo07, Stei03}. 

The $k$-path partition problem is solved for trees \cite{YaCh97}, cacti \cite{JiLi06}, cographs for a fixed $k$ \cite{Stei00}, threshold graphs \cite{Stei00}, bipartite permutation graphs \cite{Stei03}. 
The $3$-path partition problem is NP-complete on general graphs \cite{GaJo79} and comparability graphs \cite{Stei03}. The $k$-path partition problem is NP-complete on the class of cographs \cite{Stei00} and chordal bipartite graphs \cite{Stei03} if $k$ is part of the input. The $k$-path partition problem as well as the induced $k$-path partition problem remains NP-complete on bipartite graphs of maximum degree three \cite{MoTo07}. The $k$-path partition problem remains NP-complete for a graph class defined by finitely many minimal forbidden induced subgraphs \cite{Korp18}. Korpelainen \cite{Korp18} has nicely analyzed the differences between the $k$-path cover problem and induced $k$-path cover problem. The author also provides a good account of  literature survey on $k$-path cover problems and their related topics \cite{Korp18}.

\subsection{The $S$-path cover problem}
\label{subsec:S-path_cover_prob}

Here is another generalization of the path cover problem. For a set $S$ of vertices in a graph $G = (V, E)$, an \emph{$S$-path cover} is a path cover $P$ in which every vertex of $S$ is an endpoint of a path in $P$. The \emph{$S$-path cover problem} is to find an $S$-path cover of minimum cardinality in $G$. Some authors \cite{AsNi07b, AsNi10, Bake13, LiWu17} call it \emph{$k$-fixed-endpoint path cover problem} (or \emph{$kPC$ problem}) where $k = |S|$. Notice that, when $S$ is empty, the $S$-path cover problem coincides with the classical path cover problem \cite{YeCh98}. 

It is shown that the $1$-Fixed-Endpoint path cover problem is polynomial on interval graphs \cite{LiWu17} and the $k$-fixed-endpoint path cover problem is polynomial on cographs \cite{AsNi07b}, proper interval graphs \cite{AsNi10}. The $S$-path partition number is polynomially solvable for bipartite distance-hereditary graph \cite{YeCh98}. The PhD thesis by Baker \cite{Bake13} is a good source of information on this problem.

\subsection{The $k$-path vertex cover problem}
\label{subsec:k-path_ver_cover_prob}

There is another concept in the literature called \emph{$k$-path vertex cover} (which is denoted by $k$-(All-) Path Cover) \cite{FuNu14}. A set $C$ of vertices is a \emph{$k$-path vertex cover} if every $k$-path on $k$ vertices contains at least one vertex from $C$ \cite{LiZu18}. While the $k$-path cover problem is an extension of edge cover problem, the $k$-path vertex cover problem is an extension of vertex cover problem. Brause and Krivo\v{s}-Bellu\v{s} \cite{BrKr17} have studied the relationship between the $k$-path vertex cover problem and the $k$-path partition problem. Both these problems are completely different from the computational complexity point of view. The $k$-path vertex cover problem is NP-complete for cubic planar graphs of girth 3 for $k = 3$ \cite{TuYa13}.  Funke, Nusser and  Storandt \cite{FuNu14} have provided a list of applications of this problem on different domains.
 
\section{A discussion on the Gallai-Milgram theorem}
\label{sec:dis_Gallai_Milgram_thm}
One of the earliest contributions to the path partition problem is the Gallai-Milgram theorem. While surveying the Gallai-Milgram theorem, Berge  states that the Gallai-Milgram theorem has been independently extended by Las Vergnas and Linial \cite{Berg83, Hart88, Hart06, Lini78}. The Gallai-Milgram theorem illustrates the differences among the three path-types such as path, induced path and isometric path. The Gallai-Milgram theorem is true only for path covers and path partitions in undirected graphs. Let $\alpha(G)$ denote the independence number of $G$.
\begin{theorem}[Gallai-Milgram]
	\label{TGM1}
	Let $P$ be a minimal path partition of an undirected graph $G$. Then there is an independent set $S$ in $G$ such that $S$ contains exactly one vertex from each path in $P$. The vertex set $V$ of $G$ can be partitioned by at most $\alpha(G)$ paths. 
\end{theorem}

\begin{theorem}[Gallai-Milgram]
	\label{TGM2}
	Let $G$ be an undirected graph. Then $\pi_{c}(G)$ $\leq$ $\pi_{p}(G)$ $\leq$ $\alpha(G)$.
\end{theorem}
\begin{figure}[ht!]
	\begin{center}
		\scalebox{0.8}{\includegraphics{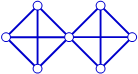}}
	\end{center}
	\caption{The Gallai-Milgram theorem fails for induced paths and isometric paths in undirected graphs}
	\label{FCExGM}
\end{figure} 
The Figure \ref{FCExGM} demonstrates how the Gallai-Milgram theorem is NOT true for induced path cover and isometric path cover. Let us make some interesting observations on the Gallai-Milgram theorem. While the Gallai-Milgram theorem fails for the induced path cover and the isometric path cover on chordal graphs, it is true on bipartite graphs.

\begin{theorem}
	\label{TGMonCho}
	Let $G$ be an undirected chordal graph without any pendant vertices. Then
	
	$\alpha(G) \leq \rho_{c}(G) \leq \rho_{p}(G)$;
	
	$\alpha(G) \leq ip_{c}(G) \leq ip_{p}(G)$;	
\end{theorem}
\proof
In chordal graphs, $\rho_{c}(G) = ip_{c}(G)$ and $\rho_{p}(G) = ip_{p}(G)$. It is enough to prove that $\alpha(G) \leq \rho_{c}(G)$. The inequality can be proved by induction on the number $q$ of maximal cliques of a chordal graph $G$. The base case is a single clique. Applying perfect elimination ordering, a few simplicial vertices (a minimum of two vertices) can be removed from $G$ in such a way that the remaining graph $G'$ is a chordal graph with $q-1$ maximal cliques. Since there are no pendant vertices of $G$, it is necessary to remove a minimum of two vertices from $G$ to build $G'$. Since $\alpha(G') \leq \rho_{c}(G')$ by induction assumption, it follows that $\alpha(G) \leq \rho_{c}(G)$.
\qed

Notice that the inequalities in Theorem \ref{TGMonCho} are not true for trees. The following result is straightforward and does not require any proof:

\begin{theorem}
	\label{TGMonBip}
	Let $G$ be an undirected bipartite graph. Then
	
	$\rho_{c}(G) \leq \rho_{p}(G) \leq \alpha(G)$;
	
	$ip_{c}(G) \leq ip_{p}(G) \leq \alpha(G)$
\end{theorem}
In undirected graphs, the Gallai-Milgram theorem illustrates how the behavior of isometric and induced paths changes from the dense graphs (such as chordal graphs without pendant vertices) to sparse graphs (such as bipartite graphs). 

Following Las Vergnas  and Linial, Berge takes the Gallai-Milgram theorem to the next level which is the famous Berge path partition conjecture \cite{Berg83, Hart88, Hart06} and weak path partition conjecture \cite{Lini81}. Hartman \cite{Hart06} has provided a survey on the Berge path partition conjectures. 

As we have discussed the behavior of the Gallai-Milgram theorem on bipartite graphs and chordal graphs, a similar discussion on the Berge path partition conjecture will be interesting. It will also be useful to study the behavior of the Berge path partition conjecture for different path-types such as paths, induced paths, and isometric paths on special classes of undirected graphs.

\section{The three cover problems}
\label{sec:three_cover_prob}
In this section, we shall discuss three problems: the path cover problem, the induced path cover problem and the isometric path cover problem.
\subsection{The path cover problem}
\label{subsec:path_cover_prob}
The Hamiltonian path problem, and hence the path cover problem is NP-complete for planar graphs, bipartite graphs, chordal graphs, chordal bipartite graphs, and strongly chordal graphs \cite{YeCh98}. Lin, Olariu and  Pruesse \cite{LiOl95} have presented an optimal algorithm for determining a minimum path cover for cographs.
Unfortunately, there is no record of any significant research work on the path cover problem. But there is a variant of path cover problem which is called the \emph{identifying path cover problem}. A path cover $\widetilde{P}$ is an \emph{identifying path cover} if for each pair $u$, $v$ of vertices, there is a path $P$ of $\widetilde{P}$ such that exactly either $u$ or $v$ belongs to $P$. Foucaud and Kov\v se \cite{FoKo13} have given lower and upper bounds on the minimum size of an identifying path cover for general graphs and discuss the tightness of the bounds. They have also illustrated that any connected graph $G$ has an identifying path cover of size at most $\lceil 2|V|/3 \rceil + 5$.
\subsection{The induced path cover problem}
\label{subsec:induced_path_cover_prob}
The induced path cover problem is NP-complete for general graphs \cite{GaJo79}. There is no record of any significant research work on the induced path cover problem too.
\subsection{The isometric path cover problem}
\label{subsec:iso_path_cover_prob}
There may not be vast difference between paths and induced paths. But there is distinct difference between induced paths and isometric paths. The following problems illustrate the difference: 
\begin{enumerate}
	\item[(1)] The problem of finding an induced path between a vertex $s$ and a vertex $t$ containing a vertex $v$ remains NP-complete in bipartite graphs \cite{Mezz10}. However, the problem of finding an isometric path between a vertex $s$ and a vertex $t$ containing a vertex $v$ is polynomially solvable for general graphs. 
	\item[(2)] Given a positive integer $k$, the problem of finding an induced path of length at least $k$ is NP-complete \cite{GaJo79}. However, the problem of finding an isometric path of length at least $k$ is a trivial problem.  
	\item[(3)] The problem of deciding, given a graph $G$ and two vertices $s$ and $t$, whether there exists an induced path of given parity between $s$ and $t$ in $G$ is NP-complete \cite{KaNi12}. But the corresponding isometric problem is polynomially solvable. 
\end{enumerate}

Fitzpatrick \cite{Fitz97} has introduced the concept of the isometric path cover problem in her PhD thesis \cite{PaCh05a}. The isometric path cover number is computed for tress, cycles, complete bipartite graphs, the Cartesian product of paths (including hypercubes) under some restricted cases \cite{FiFi01, Fitz97, Fitz99, FiNo01}. Fisher and Fitzpatrick \cite{FiFi01, FiJa06} have derived a lower bound that $ip_{c}(G) \geq \lceil |V|/({\rm diam}(G)+1) \rceil$ and have shown that the isometric path cover number of an $r \times r$ grid is $\lceil 2r/3\rceil$ and the isometric path cover number of a tree with $\ell$ leaves is $\lceil \ell /2\rceil$ . Fitzpatrick et al. \cite{FiNo01} have shown that the isometric path number of hypercube $Q_r$ is at least $2r/(r+1)$. In addition, they have also shown that $ip_{c}(Q_r) = 2^{r-{\log}_{2}(r+1)}$ when $r+1$ is a power of 2. Pan and Chang have given a linear-time algorithm to solve the isometric path cover problem on block graphs \cite{PaCh05a}, complete $r$-partite graphs and Cartesian products of 2 or 3 complete graphs \cite{PaCh06}. 
 
The strong geodetic problem is a variation of the geodetic problem which is to find a minimum set $S$ of vertices of a graph $G$ such that each vertex of $G$ lies on some isometric path between a pair of vertices from $S$ \cite{HaLo93}. A \emph{strong geodetic set} of a graph $G$ is a set $S$ of vertices such that each vertex of $G$ lies on a fixed isometric path between a pair of vertices from $S$ \cite{MaKl17}. For a graph $G$, its \emph{strong geodetic number} $sg(G)$ is the cardinality of a minimum strong geodetic set of $G$. The relationship between the isometric path cover number and the strong geodetic number is given below: 
\begin{theorem}[\cite{MaKl17}] Let $G$ be a graph.
\begin{enumerate}
	\item[(1)] When $G$ is an undirected graph, $\frac{1+\sqrt{(8\times ip_{c}(G)+1)}}{2}\leq sg(G) \leq 2\times ip_{c}(G)$.
	\item[(2)]  When $G$ is a tree or block graph, $ip_{c}(G) = \lceil sg(G)/2\rceil$.  
\end{enumerate}
\end{theorem}

Even though the isometric path cover problem is a fundamental concept in graph theory, to our knowledge, there are no literatures on the computational complexity status of the isometric path cover problem for general graphs. 

\section{The three partition problems}
\label{sec:three_part_prob}
In this section, we discuss the other set of three problems: the path partition problem, the induced path partition problem and the isometric path partition problem.
\subsection{The path partition problem}
\label{subsec:path_part_prob}
The path partition problem was introduced by Ore in 1961 \cite{Ore61}. The Hamiltonian path problem is NP-complete and hence the path partition problem is NP-complete for planar graphs, bipartite graphs, chordal graphs, chordal bipartite graphs, and strongly chordal graphs \cite{YeCh98}. Nakano, Olariu and Zomaya \cite{NaOl03} call this problem as notoriously difficult and have reported an NC-algorithm to solve the path partition problem for cographs. Dinh \cite{Dinh97} has studied the path partition problem for  graphs with toughness greater than or equal to $1$.   Asdre,  Nikolopoulos and Papadopoulos \cite{AsNi07a} have presented an optimal parallel algorithm for the path partition problem on $P_4$-sparse graphs. Pan and Chang \cite{PaCh05b} have designed a linear-time algorithm for the path-partition problem in graphs whose blocks are complete graphs, cycles or complete bipartite graphs. The path-partition problem is polynomially solvable for forests \cite{Skup74}, interval graphs \cite{ArRa90, Chan08, HuCh11}, circular arc graphs \cite{HuCh06}, distance-hereditary graphs \cite{HuCh07}, bipartite permutation graphs \cite{SriSu93}, block graphs \cite{Wong99, YaCh94}, cographs \cite{ChKu96}, bipartite distance-hereditary graphs \cite{YeCh98}, co-comparability graphs \cite{CoDa13, KoMo14}. Given a connected graph $G$, there is a spanning tree $T$ of $G$ such that $\pi_{p}(G) = \pi_{p}(T)$ \cite{BoCh74}. The relationship between the path partition number and $L(2, 1)$-labeling number of graphs have been studied by several authors \cite{ChKu96, LuZh13}.

\subsubsection{Graffiti.pc conjecture}
\label{subsubsec:graffiti_pc_conj}
The number of edges in a maximum matching of a graph $G$ is called the matching number of $G$, denoted by $\alpha'(G)$. The total domination number of $G$, denoted by $\gamma_{t}(G)$, is a minimum set $S$ of vertices such that every vertex of $G$ is adjacent to some vertex in $S$. The Graffiti.pc conjectures are well-known for total domination number \cite{DeLi07}. Some of the conjectures are related to path partition number. DeLaVi\~ na et al. \cite{DeLi07} have positively answered to some Graffiti.pc conjectures which are related to path partition numbers.
\begin{theorem}[\cite{DeLi07}]
	\label{TDeLi07}
	For every graph $G$ with no isolated vertex, $\gamma_{t}(G) \leq \alpha'(G) + \pi_{p}(G)$. When $G$ is a connected $3$-regular graph, $\gamma_{t}(G) \geq 2\pi_{p}(G)$.
\end{theorem}
Henning and Wash \cite{HeWa17} have improved Graffiti.pc conjectures and shown that the bounds are tight.
\begin{theorem}[\cite{HeWa17}]
	\label{THeWa17}
	If $G$ is a graph of order $n$, then $\alpha'(G) + (1/2)\pi_{p}(G) \geq (n/2)$. Moreover, when $\delta(G) \geq 3$, then $\gamma_{t}(G) \leq \alpha'(G) + (1/2)(\pi_{p}(G) - 1)$.
\end{theorem}

\subsubsection{Some bounds on the path partition number $\pi_{p}(G)$}
\label{subsubsec:bound_path_part_num}
While introducing the path partition problem, Ore \cite{Ore61} has produced the following bound: 
\begin{theorem}[\cite{Ore61}]
	\label{TOre61}
	Given a graph $G$ of order $n$, the path partition number $\pi_{p}$ satisfies $\pi_{p}(G) \leq n - \sigma_{2}(G)$, where $\sigma_{2}(G)$ is the minimum sum of degrees of two nonadjacent vertices.
\end{theorem}
According to Magnant and Martin \cite{MaMa09}, the above bound by Ore is sharp. There is another strong bound on the path partition number due to Hartman \cite{Hart88}. The Gallai-Milgram theorem with respect to the path partition number which is $\pi_{p}(G) \leq \alpha(G)$ is improved by Hartman \cite{Hart88} under certain conditions:

\begin{theorem}[\cite{Hart88}]
	\label{THart88}
	Let $G$ be a graph with connectivity $k$. If $\alpha(G) > k$, then the path partition number satisfies $\pi_{p}(G) \leq \alpha(G)-k$.
\end{theorem}
Magnant and Martin \cite{MaMa09} have posed the following conjecture:
\begin{conjecture}[\cite{MaMa09}]
	\label{CMaMa09}
	Given a $k$-regular graph $G$ of order $n$, the path partition number satisfies $\pi_{p}(G) \leq n/(k + 1)$.
\end{conjecture}
The bound in the Conjecture \ref{CMaMa09} is sharp \cite{MaMa09}. In support of their conjecture, they have provided the following result:
\begin{theorem}[\cite{MaMa09}]
	\label{TMaMa09}
	Given a $k$-regular graph $G$ of order $n$ with $0 \leq k \leq 5$, the path partition number satisfies $\pi_{p}(G) \leq n/(k + 1)$.
\end{theorem}
Magnant, Wang and Yuan \cite{MaWa16} have recently revised the Conjecture \ref{CMaMa09} as follows:
\begin{conjecture}[\cite{MaWa16}]
	\label{CMaWa16}
	Given positive integers $\delta$ and $\varDelta$ with $\delta \leq \varDelta$, if $G$ is a graph of order $n$ with $\delta \leq \delta(G) \leq \varDelta(G) \leq \varDelta$, then
	$\pi_{p}(G) \leq \max\{n/(\delta+1), (\varDelta-\delta)n/(\varDelta+\delta) \}$.
\end{conjecture}
They have also proved this conjecture when $\delta = 1, 2$. They have also provided another bound as follows:
\begin{theorem}[\cite{MaWa16}]
	\label{TMaWa16}
	Suppose $G$ is a graph of order $n$ with $\varDelta(G)=\varDelta < 2\delta = 2\delta(G)$, then 
	$\pi_{p}(G) \leq ((\varDelta-2)n/(2(\varDelta+\delta-4))$. 
\end{theorem}
While pointing out that Conjecture \ref{CMaMa09} is true even for $k \geq (n-1)/2$, Han \cite{Han18} has provided a tight bound for bipartite graphs as follows:
\begin{theorem}[\cite{Han18}]
	\label{THan18}
	Given a bipartite $k$-regular graph $G$ of order $n$, the path partition number satisfies $\pi_{p}(G) \leq n/2k$. 
\end{theorem}
For $3$-regular graphs, Reed \cite{Reed96} has provided a bound which is sharper than the bound in Conjecture \ref{CMaMa09}.
\begin{theorem}[\cite{Reed96}]
	\label{TReed96}
	Given a connected $3$-regular graph $G$ of order $n$, the path partition number satisfies $\pi_{p}(G) \leq n/9$. 
\end{theorem}
Reed \cite{Reed96} also demonstrates that these results are tight. Furthermore, Reed has conjectured in the same paper that 
\begin{conjecture}[\cite{Reed96}]
	\label{CReed96}
	Given a $2$-connected $3$-regular graph $G$ of order $n$, the path partition number satisfies $\pi_{p}(G) \leq \lceil n/10\rceil$. 
\end{conjecture}
This conjecture has been recently settled by Yu \cite{Yu17}. In addition, Yu has given an interesting lower bound for $\pi_{p}(G)$:
\begin{theorem}[\cite{Yu17}]
There are infinitely many $2$-connected $3$-regular $n$-vertex graphs whose path cover numbers are at least $n/14$.
 
\end{theorem}
\subsection{The induced path partition problem}
\label{subsec:indu_path_part_prob}
Chartrand et al. \cite{Char94} have introduced the concept of induced path partition problem and have produced the induced-path number of complete bipartite graphs, complete binary trees, $2$-dimensional meshes, butterflies, and general trees. They have shown that $\rho_{p}(P_m \times P_n) = 2$. In addition, they have also derived upper bounds for hypercubes and conjectured that $\rho_{p}(Q_{r}) \leq r$ for the $r$-dimensional hypercube $Q_r$ with $r \geq 2$ \cite{Char94}. Alsardary \cite{Alsa97} has answered to the conjecture that $\rho_{p}(Q_r) \leq 16$ \cite{PaCh07}. The induced path partition problem is concluded to be NP-complete by showing that it is NP-complete to determine whether the vertex set of a graph can be partitioned into two induced paths \cite{LeLe03}. Pan and Chang \cite{PaCh07} have presented a linear-time algorithm for the induced path partition problem on graphs whose blocks are complete graphs, cycles or complete bipartite graphs. Broere, Jonck and Domke \cite{BrJo06} have determined $\rho_{p}(G)$ for $G = K_m \times K_n$ and have shown that $\rho_{p}(C_m \times C_n) < 3$. 

Broere et al. \cite{BrJo05} have studied the induced path partition problem on the complements of some graphs such as cycles, paths, product of complete graphs, and product of cycles. They have also produced an elegant Nordhaus-Gaddum type relation in terms of induced path partition number $\rho_{p}$ as follows:

\begin{theorem}[\cite{BrJo05}]
	\label{TBrJo05}
	Given an undirected graph on n vertices, 
	\begin{enumerate}
		\item[(1)] $\sqrt{n} \leq \rho_{p}(G) + \rho_{p}(\overline{G}) \leq \lceil3n/2\rceil$ for arbitrary graph $G$
		\item[(2)] $1 + \lceil n/4 \rceil \leq \rho_{p}(G) + \rho_{p}(\overline{G}) \leq \lceil3n/2\rceil$ for bipartite graph $G$.
		\item[(3)] $|\ell| \leq \rho_{p}(T) + \rho_{p}(\overline{T}) \leq \lceil3(n-1)/2\rceil$ for trees $T$ where $\ell$ is the \# of leaves of $T$. 
	\end{enumerate}
\end{theorem}

Several authors \cite{AIM08, Hogb2010, PeFo18, Row11, TaFa14} have studied the relationship between the zero-forcing number $z(G)$ and the induced path partition number $\rho_{p}(G)$. Hogben \cite{Hogb2010} have shown that $\rho_{p}(G) \leq z(G)$ for general graphs. It is shown that $\rho_{p}(G) = z(G)$ when $G$ is a block-cycle graph, unicyclic graph, double path, and series of double paths \cite{TaFa14}, a tree \cite{AIM08} and a cactus \cite{Row11}.

Catral et al. \cite{CaCe12} have studied the induced path partition problem for edge subdivision graphs. Let $e = uv$ be an edge of $G$. An edge subdivision graph $G_e$ is obtained from $G$ by inserting a new vertex $w$ into $V$, deleting the edge $e$ and inserting edges $uw$ and $wv$. 

Let us recall that $M(G)$ denotes the maximum nullity of $G$. Johnson and Duarte \cite{JoDu99} have shown that $M(G) = \rho_{p}(G)$ when $G$ is a tree. Barioli,  Fallat and Hogben \cite{BaFa05} have proved that $M(G) = \rho_{p}(G)$ or $M(G) = \rho_{p}(G)-1$ for unicyclic graphs. Sinkovic \cite{Sink10} have derived that $M(G) \leq \rho_{p}(G)$ for any outerplanar graph $G$ and $M(G) = \rho_{p}(G)$ for any partial $2$-path $G$. Barioli, Fallat and Hogben \cite{BaFa04} have studied the relationship between minimum rank and the induced path partition number and have derived that $\varDelta(G) \leq \rho_{p}(G)$ where $\varDelta(G) = \max \{p-q /$ there are $q$ vertices of $G$ whose deletion leaves $p$ induced paths$\}$. 

\subsection{The isometric path partition problem}
\label{subsec:iso_path_part_prob}
It is quite a big surprise that there is no literature on the isometric path partition problem including the computational complexity status of isometric path partition problem for general graphs. 

\section{Conclusion}
\label{sec:conclusion}
This article is meant for young researchers in graph theory. It identifies the gaps and open space in ``covering and partitioning by path-types'' as potential research topics. 
Among all the six problems, the most attractive and interesting problem is path partition problem. The literatures on the path partition problem are larger than those of the other five problems. The reason may be due to the Gallai-Milgram theorem and Berge’s path partition conjecture. 
There are a very few literatures in path cover problem, induced path problem and isometric partition problem. 
Between the covering and partitioning problems, it is observed that researchers are more active on path-type partition problems than path-type cover problems. 
 
In addition to the six problems, we have also considered several problems which are related to these six problems. We have also discussed the correlation among these six problems and their related problems. We have covered the literatures for the last 30 years. 

Chang and his teammates have carried out significant research on paths \cite{ChKu96, PaCh05b, YaCh94, YaCh97, YeCh98}, induced paths \cite{PaCh07} and isometric paths \cite{PaCh05a, PaCh06}. Young researchers are encouraged to read these papers to understand how Chang and his teammates apply different techniques to solve different path-type problems.  

We have discussed only vertex version of the covering and partitioning problems. The edge version of these problems is equally popular among the researchers. In the same way, covering and partitioning by cycles or trees are other well-known subtopics of research in graph theory.
\section*{Acknowledgment}
The author expresses his deep gratitude to Prof Sandi Klav\v{z}ar for his valuable contribution and encouragement. The author also thanks Prof Michael Henning and Prof Gexin Yu for their comments.

\providecommand{\bysame}{\leavevmode\hbox to3em{\hrulefill}\thinspace}
\providecommand{\MR}{\relax\ifhmode\unskip\space\fi MR }
\providecommand{\MRhref}[2]{%
	\href{http://www.ams.org/mathscinet-getitem?mr=#1}{#2}
}
\providecommand{\href}[2]{#2}

\end{document}